\documentclass[12pt,twoside]{article}
\usepackage{amssymb}
\usepackage{amsmath,amscd}
\usepackage{pst-node,pstricks,multido,pst-plot,pst-text,pst-3d}%
\usepackage{graphicx}
\usepackage{amsfonts}

\font\elevensym=msbm10 at 11pt
\font\eightsym=msbm8
\font\sixsym=msbm6
\newfam\symfam
\textfont\symfam=\elevensym
\scriptfont\symfam=\eightsym
\scriptscriptfont\symfam=\sixsym
\def\sym{\fam\symfam\elevensym}
\def\bbbr{{\sym R}}

\def\bbbq{{\sym Q}}
\def\bbbn{{\sym N}}
\def\bbbc{{\sym C}}
\def\bbbz{{\sym Z}}
\def\bbbb{{\sym B}}
\def\qed{$\Box$}
\newtheorem{remarks}{Remarks}
\newtheorem{example}{Example}[section]
\newtheorem{rem}[example]{Remark}
\newtheorem{thm}[example]{Theorem}
\newtheorem{cor}[example]{Corollary}
\newtheorem{defn}[example]{Definition}
\newtheorem{prop}[example]{Proposition}

\newtheorem{lem}[example]{Lemme}

\newcommand{\splus}{\framebox(6.5,6.5){$+$}}
\newcommand{\sdot}{\framebox(6.5,6.5){$.$}}

\newcommand{\scriptinf}{\framebox(4.7,9){\scriptsize{$\uparrow$}}}
\newcommand{\scriptplus}{\framebox(4.7,4.7){\scriptsize{$+$}}}
\newcommand{\scriptdot}{\framebox(4.7,4.7){\scriptsize{$.$}}}
\newcommand{\scriptstar}{\framebox(4.7,4.7){\scriptsize{$*$}}}
\def\boite{\hbox{\rlap{$\sqcap$}$\sqcup$}}
\def\OV#1{\overline{#1}}

\newcommand{\scriptshuffl}{\framebox(12,6){\scriptsize{$\shuffl$}}}

\newcommand{\shuffl}{
\kern 1pt
\rule{.3pt}{4pt}\rule{4.5pt}{.3pt}\rule{.3pt}{4pt}\rule{4.5pt}{.3pt}\rule{.3pt}{4pt}
\kern 1pt
}
\begin{document}

\title{ Direct and dual laws for automata with multiplicities}
\author{G. Duchamp, M. Flouret, \'E. Laugerotte and J-G. Luque\\
LIFAR, Facult\'e des Sciences et des Techniques,\\ 76821
Mont-Saint-Aignan CEDEX France.}

%
\maketitle
\begin{abstract}
We present here theoretical results coming from the implementation
of the package called AMULT (automata with multiplicities in several
noncommutative variables). We show that classical formulas are
``almost every time'' optimal, characterize the dual laws preserving
rationality and also relators that are compatible with these laws.
\end{abstract}

\section{Introduction}
Noncommutative formal series (i.e. functions on the free monoid, with values
in
a - commutative or not - semiring) encode an infinity of data.
Rational series can be represented by linear recurrences, corresponding to
automata with multiplicities, and therefore they can be generated by finite
state processes.
Literature can be found on these ``weighted automata'' and their
theoretical and practical (e.g. \cite{KS}, \cite{SS}, \cite{HK}, \cite{CK}, \cite{MPR})
applications
(recently one of us solved
a conjecture in operator theory using these tools \cite{DR}).
The theory was founded by
Sch\"utzenberger in 1961 \cite{Sc1} where the link between recognizable and
rational series is showed (see also \cite{Sc4}), extending to rings (and to
semirings \cite{BR}) Kleene's result for languages \cite{Kl}
(corresponding to boolean coefficients). In 1974, for the case of fields,
Fliess \cite{Fl1} extended the proof of the equivalence of minimal linear
representations, using Hankel matrices. All these results allow us to
construct an algorithmic processing for this series and their associated
operations. In fact, classical constructions of
language theory have multiplicity analogues which can be used in every
domain
where linear recurrences between words are handled.
All these operations can be found in the package over automata with
multiplicities (called AMULT). This package is a
component of the environment SEA (Symbolic Environment for Automata) under
development at the University of Rouen.

The structure of this paper is the following: In section $3$ (the first section after introductory paragraphs), we recall the classical
construction for simple rational laws ($+,.,*,\times$) and make some remarks concerning in particular the non-commutative case.
The compositions are based on polynomial formulas which has an important consequence on composition of automata choosen
"at random". In fact, this first result says that the classical formulas are "almost everywhere" optimal (which is clear from
experimental tests at random).

In section $4$, we show that the three laws known to preserve rationality ( Hadamard, shuffle and infiltration products) are of the
same nature: they arise by dualizing alphabetic morphisms. Moreover, they are, up to a deformation, the only ones of this kind,
which of course, shows immediately in the implemented formulas.

Section $5$ is devoted to study the compatibility with relators. It was well known that, when coefficients are taken in a ring of
characteristic $0$, the only relators compatible with the shuffle were partial commutations 
(\cite{DK1}). Here, we show that a similar result
holds (up to the supplementary possibility of letters erasure) when $K$ is a semiring which is
not a ring. This implies the known case as a corollary. To end with, we give examples of some strange relators in
characteristic $2$.
\section{Preamble}
Let $K \langle \langle A\rangle \rangle$ be the set of noncommutative formal
series with $A$ a finite alphabet and $K$ a semiring (commutative or not). A series denoted
$ S = \sum_{w \in A^*} \langle S | w \rangle w $ is recognizable iff
there exists a row vector
$\lambda \in K^{1 \times n}$, a morphism of monoids $\mu : A^* \rightarrow
K^{n \times n}$
and a column vector $\gamma \in K^{n \times 1}$, such that for all $w \in
A^*$,
one has $\langle S |w \rangle = \lambda \mu (w) \gamma$.
Throughout the paper, we will denote by
$S:(\lambda, \mu, \gamma)$ this property and say that $(\lambda, \mu,
\gamma)$
is a linear representation of $S$, or an automaton with behaviour $S$. The integer $n$ is called
the
{\it dimension} of the linear representation $(\lambda, \mu, \gamma)$ \cite{Fl1}.

Let $K^{\mbox{\scriptsize rat}}\langle \langle A \rangle \rangle$ be the set
of rational noncommutative formal series, that is the set generated from the
letters and the laws ``$.$'' (concatenation or Cauchy product), $*$ (star
operation, partially defined), $\times$ (external product) and $+$ (union or
sum). The preceding four laws are called simple rational laws. The following important theorem for series \cite{Sc1} is the
analogue
of Kleene's theorem for languages (and in fact implies it).
\begin{thm}[Sch\"utzenberger, 1961]
A formal series is recognizable if and only if it is rational.
\end{thm}
Notice that, in the boolean case, $\times$ (the external product) is trivial, but it permits to take for granded that $L=\emptyset$
and then $\emptyset^*=1$ are rational (see \cite{Kl,HH}).

A reduced automaton $(\lambda, \mu, \gamma)$ is an automaton of minimal dimension among all the
automata with behaviour $S$\footnote
{Existence is assumed by definition, unicity is proved in case $K$ is
$\bbbb$ (for deterministic automata) or a
(commutative or not) field \cite{FL} but is problematic in general.}. This
minimum is called the {\it rank} of the series $S$ \cite{Sc1}. In case $K$
is
a field, the rank of $S$ is the dimension of the linear span
of the shifts of $S$ (see Sect. \ref{sec3}). It is the smallest number of
nodes of an automaton with behaviour $S$. Here,
minimization (up to an equivalence) is possible \cite{Sc1} (see also
\cite{BR}). An explicit
algorithm is given in full details in \cite{FL} (notice that this algorithm
is valid as well for noncommutative multiplicities) as well as the
construction
of intertwining matrices.

Again, the specialisation of $K$ to the boolean semiring $\bbbb$ yields
to
the case of classical finite state automata.

\section{Constructing usual laws}\label{sec3}
\subsection{Operations on linear representations}
We expound here universal formulas for constructing linear representations.
They can be applied to any
semiring $K$. For two representations of ranks $n$ and $m$,
it will be provided a representation of rank $r(n,m)$.
Let us recall some classical facts.
Classical operations on series are sum, external product and star (unary and
partially defined). By definition, the sum of two series $R$ and $S$ is
$$
R+S = \sum_{w \in A^*}
\bigl(\langle R|w\rangle + \langle S|w\rangle \bigr) w,
$$
their concatenation (or Cauchy product)
$$
R.S = \sum _{w \in A^*} \left( \sum _{uv=w}\langle R|u\rangle \langle
S|v\rangle \right) w,
$$
and the star of a series $S$
$$
S^* = \sum_{n \geq 0} S^n= 1 +S S^*
$$
if its constant term is zero (such a series is said to be proper). The
preceding
operations have polynomial counterparts in terms of linear
representations. We gather them in the following proposition.
\begin{prop} \label{p1}
Let $R : {\cal A}_{r}=(\lambda^r,\mu^r,\gamma^r)$ (resp.
$S:{\cal A}_{s}=(\lambda^s,\mu^s,
\gamma^s)$) of rank $n$ (resp. $m$). The linear representations of the sum,
the
concatenation and the star are respectively \\
$R+S:$
\begin{eqnarray}\label{sum}
{\cal A}_{r}\mbox{ } \splus \mbox{ } {\cal A}_{s} =\left( \left(
\begin{array}{cc} \lambda^r& \lambda^s
\end{array} \right),
\left(
\begin{array}{c|c}
\mu^r(a)&0_{n\times m} \\ \hline 0_{m\times n} &\mu^s(a)
\end{array}
\right)_{a \in A},
\left(
\begin{array}{c}
\gamma^r \\ \gamma^s
\end{array}
\right)
\right) \enspace ,
\end{eqnarray}
$R.S:$
\begin{eqnarray}\label{cau}
{\cal A}_r \mbox{ } \sdot \mbox{ } {\cal A}_s = \left( \left(
\begin{array}{cc}
\lambda^r& 0_{1\times m}
\end{array}
\right),
\left(
\begin{array}{c|c}
\mu^r(a) & \gamma^r\lambda^s\mu^s(a)\\
\hline
0_{m\times n} & \mu^s(a)
\end{array}
\right)_{a\in A},
\left(
\begin{array}{c}
\gamma^r\lambda^s\gamma^s\\
\gamma^s
\end{array}
\right)
\right) \enspace ,
\end{eqnarray}
\noindent If $\lambda^s \gamma^s = 0$,
$S^*:$
\begin{eqnarray}\label{star}
{\cal A}_s \mbox{ }^{\scriptstar} = \left(\left(
\begin{array}{cc}
0_{1\times m}&1
\end{array}
\right ),
\left(
\begin{array}{c|c}
\mu^s(a) +\gamma^s\lambda^s \mu^s(a)& 0_{m\times 1} \\
\hline
\lambda^s \mu^s(a) & 0
\end{array}
\right)_{a \in A},
\left(
\begin{array}{c}
\gamma^s \\ \
1
\end{array}
\right)
\right)\enspace .
\end{eqnarray}
\end{prop}
{\bf Proof}
Formula $(1)$ is straightforward.\\
To prove formula $(2)$, let
$(\lambda, \mu, \gamma):= {\cal A}_r \mbox{ } \sdot \mbox{ } {\cal A}_s $. One proves by induction that
$$
\mu(w) = \left(
\begin{array}{cc}
\mu^r(w) & \displaystyle \sum_{uw=w \atop v \neq 1} \mu^r(u) \gamma^r \lambda^s \mu^s (v) \\
0_{m \times n} & \mu^s (w)
\end{array}
\right),
$$
and then $\lambda \mu(w) \gamma =\displaystyle \sum_{uv=w} \lambda^r \mu^r (u) \gamma^r \lambda^s \mu^s(v) \gamma^s =
\sum_{uv=w} \langle R|u\rangle \langle S|v\rangle$. \\
Concerning the formula $(3)$, let $(\lambda^*, \mu^*, \gamma^*):= {\cal A}_s \mbox{
}^{\scriptstar}$. Again,
$$
\mu^*(w) = \left(
\begin{array}{cc}
* & 0_{m \times 1} \\
\sum_{n=1}^{|w|} \sum_{u_1\cdots u_n =w \atop u_i \neq 1} (\lambda_s \mu_s(u_1) \gamma_s)\cdots (\lambda_s
\mu_s(u_{n-1}) \gamma_s) (\lambda_s \mu_s(u_{n})) & 0
\end{array}
\right),
$$
that is \\
$
\begin{array}{rcl}
\lambda^* \mu^*(w)\gamma^* & = & \sum_{n=1}^{|w|} \sum_{u_1\cdots u_n =w \atop u_i \neq 1} (\lambda_s \mu_s(u_1)
\gamma_s)\cdots (\lambda_s \mu_s(u_{n}) \gamma_s)\\
& = &\sum_{n=1}^{|w|} \langle S^n |w| \rangle = \sum_{n\geq 0}
\langle S^n |w| \rangle = \langle S^* |w\rangle. \Box
\end{array}
$
\begin{rem}
\begin{enumerate}
\item Formulas (\ref{sum}) and (\ref{cau}) provide associative laws on
triplets.
They can be found explicitly in \cite{CK}.\\
\item Formula (\ref{star}) makes sense even when $\lambda^s \gamma^s \neq 0$
(this fact
will be used in the density result of Section \ref{secden}).\\
\item Of course if $S: (\lambda,\mu,\gamma)$ and $\alpha\in K$ then $\alpha S:= \alpha \times S:
(\alpha\lambda,\mu,\gamma)$ and $S\alpha:= S \times \alpha: (\lambda,\mu,\gamma\alpha)$.\\
\item For the sum (${\cal A}_R\mbox{ }\splus\mbox{ }{\cal A}_S$), ${\cal A}_R$ and
${\cal A}_S$ are just placed side by side.\\ \\
The product ${\cal A}_R\mbox{ }\scriptdot\mbox{ }{\cal A}_S$ has the following components

\begin{itemize}
\item[-] {\bf States:} States of ${\cal A}_R$ and ${\cal A}_S$.\\
\item[-] {\bf Inputs:} Inputs of ${\cal A}_R$.\\
\item[-] {\bf Transitions:} Transitions of ${\cal A}_R$ and ${\cal A}_S$ and,
for each letter $a$, each state $r_i$ of ${\cal A}_R$ and each state $s_j$
of ${\cal A}_S$, a new arc
$r_i\stackrel{a}{\rightarrow} s_j$
is added with the
coefficient $(\gamma_r)_i(\lambda_s\mu_s(a))_j$.\\
\item[-] {\bf Outputs:} The scalar product $\lambda_s\gamma_s$ is computed
once for all and there is an output on each $q_i$ with the coefficient
$(\gamma_r)_i\lambda_s\gamma_s$, the outputs of ${\cal A}_S$
being unchanged.
\end{itemize}

\bigskip
\noindent For ${\cal A}^{\scriptstar}$, one adds a new state $q_{n+1}$ with an input and an
output bearing coefficient $1$, every coefficient $\mu_{i,j}(a)$ is
multiplied by $(1+\gamma_i\lambda_j)$ and new transitions
$q_{n+1}\stackrel{a}{\rightarrow} q_{i}$ with coefficient
$\sum_k\lambda_k\mu_{k,i}(a)$ (i.e. the "charge"
of the state $q_i$ after reading $a$) are added.
\\ \\
In the case $K=\bbbb$, one recovers the classical boolean constructions.
\end{enumerate}
\end{rem}
\subsection{Sharpness}
Here we discuss the sharpness of the preceding constructions. Indeed,
testing
our package showed us that ``almost everytime'' the compound automata
was minimal when the data were choosen at random. The crucial point in the
proof of
Theorem \ref{thden} is the fact that certain polynomial indicators are not
trivial.
For this, we use suited examples which are gathered in the following
subsection.\\ \\
{ a) Test automata}

\vspace{3mm}
Let ${\cal B} = (S_i)_{1 \leq i \leq n}$ be a finite sequence of series
generating
a stable module and $S=\sum_{i=1}^n {\lambda}_i S_i$. It is well known that
the triplet
$$
\left(\sum_{i=1}^{n} \lambda_i e_i , \mbox{ }
\left( \left[ \mu_{i,j}(a)\right]_{1 \leq i,j \leq n}\right) _{a\in A},
\mbox{ }
\sum_{i=1}^{n} \langle S_i| 1\rangle e_i^*\right)
$$
(where $e_i:=(0,\cdots ,1,\cdots 0)$
with the entry $1$ at place $i$, $e_i^*$ the transpose of $e_i$, and
$a^{-1}S_i=\sum_{j=1}^n \bigl(\mu(a)\bigr)_{ij}S_j$
for any letter $a\in A$) is a linear representation of $S$.
Here, to each series of one variable, $S=\sum_{p \geq 0} \alpha_{p} a^p$, of
rank $n$,
over a field $K$, we associate the triplet $\tau (S)$ given by ${\cal
B}=(a^{-p}S)_{0 \leq p \leq n-1}$.
\begin{rem}\label{extalph} Of course, if $a \in A$ we consider that $S$ belongs to $K \langle
\langle A \rangle \rangle$
and this will neither affect the rank nor the following constructions.
\end{rem}
\begin{lem}
Let $\displaystyle S_{\alpha, n} =\frac{1}{(1-\alpha a)^{n}}$ and
$\displaystyle
T_n=\frac{a^{n-1 } }{1-a^n}$ be $\bbbq$-series.
\begin{enumerate}
\item The rank of $S_{\alpha,n}$, $S_{\alpha,n}+S_{\beta,m}$ ($\alpha \neq
\beta$), and
$S_{\alpha,n}. S_{\alpha,m}$ are respectively $n$, $n+m$ and $n+m$.
\item The rank of $T_n$ is $n$ and that of $T_n^*$ is $n+1$.
\end{enumerate}
\end{lem}
{\bf Proof} Straightforward. \qed\\

\noindent{b) ²Density}\label{secden}

\vspace{3mm}
The following theorem proves that, if the data are choosen ``at random'' in
bounded domains, the compound automaton is almost surely minimal. More
precisely:
\begin{thm}\label{thden}
Let $A$ be a finite alphabet and
${\cal A}_i=(\lambda_i,\mu_i,\gamma_i)$ two
automata of dimension $n_i$ ($i=1,2$),
choosen ``at random'' within bounded non trivial disks of $K$ ($K= \bbbr$ or
$\bbbc$). Then the probability
that the automaton ${\cal A}_1 \mbox{ }\splus
\mbox{ }{\cal A}_2$ (resp. ${\cal A}_1
\mbox{ }\sdot \mbox{ } {\cal A}_2$, ${\cal A}_1
\mbox{ } ^{\scriptstar}$) be minimal is $1$.
\end{thm}
{\bf Proof} The proof rests on the following lemma.
\begin{lem}
There is a polynomial mapping $P : K^{|A |\times n^2 +2n}
\rightarrow K^s$ such
that $P(\lambda, \mu, \gamma)=0$ iff $(\lambda, \mu, \gamma)$ (an automaton
of dimension $n$) is not minimal.
\end{lem}
{\bf Proof of the lemma} By a theorem of Sch\"utzenberger
\cite{Sc1}, the representation
%
%
$(\lambda, \mu, \gamma)$ is minimal iff $\lambda \mu (K \langle A
\rangle)=K^{1 \times n}$ (resp. $\mu (K \langle A \rangle)\gamma=K^{n \times
1}$). As there is a prefix (resp. suffix) subset $U \subset A^*$ (resp. $V
\subset A^*$) such that $\lambda \mu (U)$ (resp. $\mu(V)\gamma$) is a basis,
we have $U \subset A^{<n}$ (resp. $V \subset A^{<n}$). Let $A^{<n}=\{ w_1:=1,
w_2, \cdots , w_m\}$ ($m=(|A|^{n}-1 )/ (|A|-1)$), one constructs the $m \times n$
(resp. $n \times m$) matrix
$$
L= \left(
\begin{array}{c}
\lambda \mu(w_1) \\
\lambda \mu(w_2) \\
\vdots \\
\lambda \mu (w_m)
\end{array}
\right)
\left(\mbox{ resp. }
M = \left(
\begin{array}{ccc}
\mu(w_1) \gamma & \cdots & \mu(w_m)
\end{array}
\right)
\right),
$$
these matrices have polynomial entries in the data. In view of
what precedes, minimality is equivalent to the non nullity of some $n \times
n$-minor of $L$ and of $M$. Sorting these minors as a vector, one get the
desired polynomial mapping $ K^{|A |\times n^2 +2n}
\rightarrow K^s$ with $s=\left( {\displaystyle m} \atop {\displaystyle
n}\right)$.\qed
\\ \\
The other steps go as follows.
\begin{enumerate}
\item For the two first operations, let $P\mbox{ }_{\scriptplus}=
\left( {\cal A}_1 \mbox{ }\splus
\mbox{ }{\cal A}_2 \right)$, $P\mbox{ }_{\scriptdot}= P\left( {\cal A}_1
\mbox{ }\sdot
\mbox{ } {\cal A}_2 \right)$, and prove that $P\mbox{ }_{\scriptplus}$ (resp.
$P\mbox{ }_{\scriptdot}$) is
not trivial using $\tau (S_{\alpha,n})= {\cal A}_1 $ and
$\tau(S_{\beta,n})={\cal A}_2 $,
$\alpha \neq \beta$ (resp. $\tau (S_{\alpha,n})= {\cal A}_1 $ and
$\tau(S_{\alpha,m})
={\cal A}_2 $ ) extended to the alphabet $A$ in view of remark \ref{extalph}.
For the star operation, prove that $P\mbox{ }_{\scriptstar}= P({\cal A}_1
\mbox{ }^{ \scriptstar})$ is not trivial
using $\tau (T_n)= {\cal A}_1 $.
\item End of the proof: if $\phi : K^r \rightarrow K^s$ is polynomial and
not trivial, let $\nu$ be the normalized
uniform probability mesure on the product of disks, then the probability
such that $\phi(\nu) \neq 0$ is $1$ as
$\phi^{-1}\{0\}$ is closed with empty interior.\qed
\end{enumerate}

\section{Dual laws}

\subsection{Discussion}

Let $a,b \in A$, $u,v \in A^*$, and $\odot_{\epsilon,q}$ be the law defined recursively by
$$
\left\{ \begin{array}{l}
1 \odot_{\epsilon,q} 1 =1,\mbox{ } a \odot_{\epsilon,q} 1 = 1 \odot_{\epsilon,q} a = \epsilon a, \\
au \odot_{\epsilon,q} bv = \epsilon \bigl( a (u \odot_{\epsilon,q} bv) + b (au \odot_{\epsilon,q} v)\bigr) +
q \delta_{a,b} a (u \odot_{\epsilon,q} v)
\end{array}
\right.
$$
with $\delta_{a,b}$ the Kronecker delta.

One immediately checks that this law is associative iff $\epsilon \in \{ 0,1\}$. We get, here,
the well-known shuffle ($ \shuffl = \odot_{1,0}$), infiltration ($\uparrow = \odot_{1,1}$) and Hadamard ($\odot
=\odot_{0,1}$)
products (\cite{Ei}, \cite{Lo}). Then, $\odot_{1,q}$ is a continuous deformation between shuffle and infiltration.
These laws can be called ``dual laws'' as they proceed from the same template that we
now describe. We use an implementable realisation of the lexicographically ordered tensor
product. Let us recall that the tensor product of two spaces $U$ and $V$ with bases
$(u_i)_{i \in I}$ and $(v_j)_{j \in J}$ is $ U \otimes V $, with basis $(u_i
\otimes v_j)_{(i,j)\in I \times J}$, and for the sake of
computation, we impose that the set $I \times J$ be lexicographically ordered.

Let $ K\langle A \rangle \otimes K\langle A \rangle$ be the ``double'' non commutative
polynomial algebra that is the set of finite sums $P = \sum_{u,v \in A^*} \langle P | u
\otimes v\rangle u \otimes v$, the product being given by $(u_1 \otimes v_1)(u_2 \otimes v_2) = u_1
u_2 \otimes v_1 v_2 $.\\
The construction of dual laws is based on the following pattern:\\

Let $c:K \langle A \rangle \rightarrow K\langle A
\rangle \otimes K\langle A \rangle$,
if for all $w\in A^*$, the set
$
\{w : \langle u \otimes v |c(w) \rangle \neq 0 \}
$
is finite (in which case $c$ will be called {\it locally finite}), then the sum
$$
u \mbox{ }\boite_{\alpha}\mbox{ } v = \sum_{w \in A^*} \langle u \otimes v |c_{\alpha}(w) \rangle w
$$
exists and defines a (binary) law $\boite_{\alpha}$ on $K \langle A \rangle$, dual to
$c_{\alpha}$. Then, this extends to series by
$$
\langle R\mbox{ }\boite_{\alpha} \mbox{ } S| w\rangle :=\langle R\otimes S|c_{\alpha}(w)\rangle \enspace .
$$
One can show easily that the three laws $\odot$, $\shuffl$ and $\uparrow$ come from coproducts
defined on the words by

\begin{enumerate}
\item $c_{\alpha}(a_1a_2\cdots a_n)=c_{\alpha}(a_1)c_{\alpha}(a_2)\cdots c_{\alpha}(a_n)$,
\item $c_\odot(a)=a\otimes a$, $c_{\shuffl} (a)=a\otimes 1+1\otimes a$, $c_\uparrow (a)=a\otimes
1+1\otimes a+a\otimes a$,
\end{enumerate}

\noindent
and generally
$c_{\epsilon,q}(a)=\epsilon(a\otimes 1+1\otimes a)+q a\otimes a$.

\noindent The preceding computation scheme has an immediate consequence on the implementation of the laws.
\begin{prop}
Let $R:(\lambda^r, \mu^r, \gamma^r)$ and $S:(\lambda^s, \mu^s, \gamma^s)$. Then
$$
R \mbox{ }\boite_{\alpha}\mbox{ } S : (\lambda^r \otimes \lambda^s, \mu^r \otimes \mu^s \circ c_{\alpha},
\gamma^r \otimes \gamma^s) \enspace .
$$
\end{prop}
\noindent {\bf Proof}
We verify it by duality. Indeed, for $w \in A^*$,\\
$
\begin{array}{rl}
\langle R \otimes S | c_{\alpha}(w) \rangle =
& \sum_{u,v \in A^*} \langle \lambda^r \otimes \lambda^s \left( \mu^r \otimes \mu^s (u \otimes v) \right) \gamma^r \otimes
\gamma^s \times u \otimes v |c_{\alpha}(w)\rangle\\
= & \sum_{u,v \in A^*} \lambda^r \otimes \lambda^s \left( \mu^r \otimes \mu^s (u \otimes v) \right) \gamma^r \otimes \gamma^s
. \langle u \otimes v |c_{\alpha}(w) \rangle \\
= & \lambda^r \otimes \lambda^s \left( \sum_{u,v \in A^*} \mu^r \otimes \mu^s \langle u \otimes v |c_{\alpha}(w) \rangle (u
\otimes v) \right) \gamma^r \otimes \gamma^s \\
= & \lambda^r \otimes \lambda^s \left( \mu^r \otimes \mu^s \sum_{u,v \in A^*} \langle u \otimes v |c_{\alpha}(w) \rangle (u
\otimes v) \right) \gamma^r \otimes \gamma^s \\
= & \lambda^r \otimes \lambda^s \left( \mu^r \otimes \mu^s
c_{\alpha}(w) \right) \gamma^r \otimes \gamma^s .\Box
\end{array} $
\\

%
%
Let us study among laws which ones are associative.
\begin{prop}
Let $K$ be a field, and $c_{\alpha} : K \langle A \rangle \rightarrow K\langle A \rangle
\otimes K\langle A \rangle$ the alphabetic morphism defined on the letters of $A$ by
$$
c_{\alpha}(a) = \sum_{p,q \geq 0} \alpha_{p,q} a^p \otimes a^q
$$
with $c_{\alpha}(1) = 1 \otimes 1$ ($\alpha_{p,q} = \alpha_{p,q} (a)$ may vary from one letter to one another).
\begin{enumerate}
\item The morphism $c_{\alpha}$ is locally finite iff $\alpha_{0,0} =0$.
\item \label{iip9} Providing $\alpha_{0,0} = 0$, the following assertions are equivalent.
\begin{enumerate}
\item \label{2a} The law $ \boite_{\alpha}$ defined by $\langle u \mbox{ }\boite_{\alpha}\mbox{ } v | w\rangle :=
\langle u \otimes v | c_{\alpha }(w) \rangle$ ($u,v,w \in A^*$) is associative.
\item \label{2b} The coefficients $\alpha_{p,q}$ satisfy the relations $\alpha_{p,q} = 0$ for
$p \mbox{ or }q \geq 2$, $\alpha_{0,1}, \alpha_{1,0} \in \{ 0,1 \}$ and $\alpha_{0,1}
\alpha_{1,1} = \alpha_{1,0} \alpha_{1,1}$.
\end{enumerate}
\item Providing (\ref{iip9}.\ref{2b}), the element $1_{A^*}$ is a unit for $\boite_{\alpha}$ iff $\alpha_{0,1} =
\alpha_{1,0} = 1$.
\end{enumerate}
\end{prop}
\noindent
{\bf Proof}
\begin{enumerate}
\item We have $\displaystyle c_{\alpha}(a) = \alpha_{0,0} 1 \otimes 1 + \sum_{p+q \geq 1 } \alpha_{p,q} a^p \otimes a^q$,
and then for all $n \geq 0$, $\displaystyle c_{\alpha} (a^n)= \alpha_{0,0}^n 1 \otimes 1 + \sum_{p+q \geq 1 } \beta_{p,q} a^p
\otimes a^q$ for some $\beta_{p,q}$.
If $\alpha_{0,0}$ were not zero, the term $1 \otimes 1$ would appear in an infinity of words, and then $c_{\alpha}$ would not
be locally finite.\\
Conversely, if $\alpha_{0,0}(a)=0$ (for every letter), then $\displaystyle c_{\alpha}(a)= \sum_{p+q \geq 1} \alpha_{p,q}
a^p \otimes a^q$ and for all word $w = a_1 \cdots a_n \in A^*$,
$$
c_{\alpha}(w) = \sum_{\scriptstyle
p_i + q_i\geq 1 \atop \scriptstyle 1 \leq i \leq n} \left( \prod_{i=1}^{n} \alpha_{p_i, q_i}(a_i)\right) a_1^{p_1}\cdots
a_n^{p_n} \otimes
a_1^{q_1}\cdots a_n^{q_n}.
$$
As $p_i +q_i \geq 1$, we have $\displaystyle \sum_{i=1}^{n} (p_i +q_i) \geq n$, that is to say
$$ \langle c_\alpha(w),u\otimes v\rangle\Rightarrow \left\{
\begin{array}{l}
w| \leq |u|+|v|\\
Alph(w) = Alph(u) \cup Alph(v)
\end{array}\right.
$$
 where  $u:= a_1^{p_1}\cdots a_n^{p_n}$ and $v:=a_1^{q_1}\cdots a_n^{q_n}$.\\
To summarize, the set
$$
S=\{w/\langle u\otimes v|c_\alpha(w)\rangle\neq 0\}
$$
has bounded lengths and its alphabet is finite, $S$ is then finite.\\
\item First, remark that (\ref{iip9}.\ref{2a}) is equivalent to the condition
\begin{eqnarray}\label{equiv1}
(Id \otimes c_{\alpha}) \circ c_{\alpha} = (c_{\alpha} \otimes Id) \circ c_{\alpha}.
\end{eqnarray}
The law $\boite_{\alpha}$ is associative iff for all words $u_1$, $u_2$, $u_3$ $\in A^*$, we have
$$
(u_1 \boite_{\alpha} u_2) \boite_{\alpha} u_3 = u_1 \boite_{\alpha} (u_2 \boite_{\alpha} u_3)
$$
that is to say that, for all $w \in A^*$,
$$
\langle (u_1 \boite_{\alpha} u_2) \boite_{\alpha} u_3 | w \rangle = \langle u_1 \boite_{\alpha} (u_2 \boite_{\alpha} u_3) | w
\rangle \enspace .
$$
But one has \\
$
\begin{array}{ccl}
\langle (u_1 \boite_{\alpha} u_2) \boite_{\alpha} u_3 | w \rangle & = & \langle (u_1 \boite_{\alpha} u_2) \otimes u_3 |
c_{\alpha}(w) \rangle \\
& = & \langle u_1 \otimes u_2 \otimes u_3 | (c_{\alpha} \otimes Id) \circ c_{\alpha}(w) \rangle
\end{array}
$\\
and\\
$
\begin{array}{ccl}
\langle u_1 \boite_{\alpha} (u_2 \boite_{\alpha} u_3 ) | w \rangle & = & \langle u_1 \otimes (u_2 \boite_{\alpha} u_3) |
c_{\alpha}(w) \rangle\\
& = & \langle u_1 \otimes u_2 \otimes u_3 | (Id \otimes c_{\alpha} ) \circ c_{\alpha}(w) \rangle.
\end{array}
$\\
As $u_1$, $u_2$, $u_3$, $w$ are arbitrary, we get $(c_{\alpha} \otimes Id) \circ c_{\alpha} = (Id \otimes c_{\alpha} ) \circ
c_{\alpha}$.\\
To show the equivalence between (\ref{iip9}.\ref{2b}) and ($\ref{equiv1}$), suppose first that
($\ref{equiv1}$) holds. We endow $\bbbn ^k$ with the lexicographic order (reading from left to right for instance) which
is compatible with addition and will be denoted $\prec$ (here, $k=2,3$). Then, if it is not zero, $c_{\alpha} (a)$ can be written
$$
\alpha_{{\overline p}, {\overline q}}a^{\overline p} \otimes a^{\overline q} + \sum_{(p,q) \prec ({\overline p},{\overline q})}
\alpha_{p,q} a^p \otimes a^q \enspace ,
$$
$({\overline p},{\overline q})$ being the highest couple of exponents in the support. Then,
$\\
\begin{array}{lcl}
(c_{\alpha} \otimes Id) \circ c_{\alpha} (a) & = & \displaystyle \alpha_{{\overline p},{\overline q}} c_{\alpha}(a^{\overline
p}) \otimes a^{\overline q} + \sum_{(p,q) \prec ({\overline p},{\overline q})} \alpha_{p,q} c_{\alpha}(a^p) \otimes a^q \\
& = & \displaystyle \alpha_{{\overline p},{\overline q}}^{{\overline p}+1} a^{({\overline p})^2} \otimes a^{{\overline
p}{\overline q}} \otimes a^{\overline q} + \sum_{\scriptstyle (p,q,r)\prec \scriptstyle ({\overline p}^2,{\overline p}{\overline
q},{\overline q})} \beta_{p,q,r} a^p \otimes a^q \otimes a^r ,
\end{array}
$ \\
but
$\\
\begin{array}{lcl}
( Id \otimes c_{\alpha}) \circ c_{\alpha} (a) & = & \displaystyle \alpha_{{\overline p},{\overline q}} a^{\overline p} \otimes
c_{\alpha} (a^{\overline q}) + \sum_{(p,q) \prec ({\overline p},{\overline q})} \alpha_{p,q} a^p \otimes c_{\alpha}(a^q) \\
& = & \displaystyle \alpha_{{\overline p},{\overline q}}^{{\overline q}+1} a^{{\overline p}} \otimes a^{{\overline
p}{\overline q}} \otimes a^{(\overline q)^2} + \sum_{\scriptstyle (p,q,r)\prec \scriptstyle ({\overline p},{\overline p}{\overline
q},{\overline q}^2)} \beta_{p,q,r} a^p \otimes a^q \otimes a^r.
\end{array}
$\\
Necessarily, ${\overline p} = {\overline p}^2$ and ${\overline q} = {\overline q}^2$, which is only possible when ${\overline
p} \in \lbrace 0,1\rbrace$ and ${\overline q} \in \lbrace 0,1\rbrace$ and then $\alpha_{p,q}=0$ for $p \mbox{ or }q \geq 2$.
The equality now reads
$$
\begin{array}{ccc}
& \alpha_{1,0} a \otimes 1 \otimes 1 + \alpha_{0,1}^{2} 1 \otimes 1 \otimes a + \alpha_{0,1} \alpha_{1,1} a \otimes 1 \otimes
a & \\
& = & \\
& \alpha_{1,0}^2 a \otimes 1 \otimes 1 + \alpha_{0,1} 1 \otimes 1 \otimes a + \alpha_{1,0} \alpha_{1,1} a \otimes 1 \otimes a,
&
\end{array}
$$
which implies (\ref{iip9}.\ref{2b}). The converse is a straightforward computation.
\item The condition $1_{A^*}$ is a unit for $\boite_{\alpha}$ implies that, for $a \in A$, we have \\
$
\begin{array}{lcl}
1 \boite_{\alpha} a = a \boite_{\alpha} 1 = a & \Leftrightarrow & \langle 1 \boite_{\alpha} a| a\rangle = \langle a
\boite_{\alpha} 1 | a \rangle = 1 \\
& \Leftrightarrow & \langle 1 \otimes a | c_{\alpha}(a) \rangle = \langle a \otimes 1 | c_{\alpha}(a) \rangle = 1 \\
& \Leftrightarrow &
\left\{
\begin{array}{l}
\langle 1 \otimes a| \sum_{p,q \geq 0} \alpha_{p,q} a^p \otimes a^q \rangle = 1 \\
\langle a \otimes 1 | \sum_{p,q \geq 0} \alpha_{p,q} a^p \otimes a^q \rangle = 1
\end{array}
\right. \\
& \Leftrightarrow & \alpha_{0,1} = \alpha_{1,0} =1.
\end{array}
$\\
Conversely, the latter implies that, for each $w \in A^*$, $1 \boite_{\alpha} w = w
\boite_{\alpha} 1 = w $. \qed
\end{enumerate}

\begin{rem}
\begin{enumerate}
\item For just a commutative law the condition $\alpha_{p,q} = \alpha_{q,p}$ is sufficient. 
Moreover,
the condition (\ref{iip9}.\ref{2b}) implies $\alpha_{0,1}, \alpha_{1,0} \in \{ 0,1\}$.

\item If $\alpha_{11}\neq 0$, the only dual laws which are associative ones are $$
c_{\epsilon,q}(a)=\epsilon(a\otimes 1+1\otimes a)+q a\otimes a
$$ with
parameters $\epsilon\in\{0,1\}$ and $q \in K^\times$. Notice that in this case they are all 
commutative.

\item If $\alpha_{11}=0$, we get two degenerate laws (opposite between theimselves) which are not 
in the familly $(\boite_{\epsilon, q})$ with $\epsilon\in\{0,1\}$ and $q\in K$ corresponding to 
$\alpha_{10}=1$ and $\alpha_{10}=0$ (resp. $\alpha_{01}=0$ and $\alpha_{10}=1$). This laws are 
not commutative when $A\neq \emptyset$.

\end{enumerate}
\end{rem}

\subsection{Usual dual laws}

\noindent{ a) Shuffle and infiltration product ($\epsilon =1$, $q \in \{ 0,1\}$)}

\begin{prop} Let $R : (\lambda_1, \mu_1, \gamma_1)$ (resp. $S : (\lambda_2, \mu_2, \gamma_2)$) with
rank $n$ (resp. $m$).
\begin{enumerate}
\item Automata corresponding to shuffle and infiltration products are respectively\\
\begin{equation}\label{eq:shu}
\hspace{-1.4cm} R \shuffl S : (\lambda_1 \otimes \lambda_2,\left( \mu_1(a) \otimes I_2 + I_1
\otimes \mu_2(a)\right)_{a \in A} ,\gamma_1 \otimes \gamma_2 ) \enspace ,
\end{equation}
and
\begin{equation}\label{eq:inf}
\hspace{-1.4cm}R \uparrow S: (\lambda_1 \otimes \lambda_2, (\mu_1 (a) \otimes I_2 + I_1 \otimes \mu_2(a) +
\mu_1(a) \otimes \mu_2(a))_{a\in A} ,\gamma_1 \otimes \gamma_2) \enspace .
\end{equation}
\item \label{p2} The bound $nm$ is sharp in both cases.
\item The density result of theorem \ref{thden} holds.
\end{enumerate}
\end{prop}

{\bf Proof} Concerning point (\ref{p2}), an example reaching the
bound for any rank is to consider the families of series
$S_n=a^{n-1}$ and $T_n=b^{n-1}$ of rank $n$. The shuffle product
$S_n \shuffl S_m = a^{n-1} \shuffl b^{m-1}$ ($a\neq b \in A$) has a
minimal linear representation of rank $nm$. The same example is
valid for the infiltration product as, for $a \neq b$, $a^n\uparrow
b^m = a^n\shuffl b^m$. \qed
\\ \\
The proposition yields the following.
\begin{defn}
Let ${\cal A}_i=(\lambda_i,\rho_i,\gamma_i)$ with $i=1,2$ then we define ${\cal A}_1\mbox{ }\scriptshuffl \mbox{ }{\cal
A}_2$ and ${\cal A}_1\mbox{ }\scriptinf \mbox{ }{\cal A}_2$ by the formulas \ref{eq:shu} and \ref{eq:inf}.
\end{defn}
\begin{rem}
These laws are already associative at the level of automata.
\end{rem}

\noindent{ b) Hadamard product ($\epsilon =0$, $q=1$)}

\vspace{3mm}
We recall that the Hadamard product (\cite{Fl2}, \cite{Sc4}) of two series is the pointwise product of the corresponding
functions (on words).
We can use the machinery above to describe an automata for it.
\begin{prop}
Let $R : (\lambda^r, \mu^r, \gamma^r)$ (resp. $S : (\lambda^s, \mu^s, \gamma^s)$) with
rank $n$ (resp. $m$).
A representation of the Hadamard product is
$$
R \odot S :\left(\lambda^r \otimes \lambda^s,\left( \mu^r(a) \otimes \mu^s(a) \right)_{a\in A},\gamma^r \otimes \gamma^s \right) 
\enspace ,
$$
and the bound is asymptotically sharp.
\end{prop}
{\bf Proof} Let $\beta (n,m):=sup_{rank(R) = n \atop rank(S)=m}
rank(R\odot S)$. We claim that
$$
\limsup_{n,m \rightarrow + \infty} \displaystyle \frac{\beta (n,m)}{nm}=1 \enspace ,
$$
(what we mean by ``asymptotically sharp'').\\ Indeed, let us consider the Hadamard product of two series of the family
$$
S_n = \sum_{k \geq 0} a^{nk} = \frac{1}{(1-a^n)}\enspace .
$$
The rank of $S_n$ is $n$, and
$$
\begin{array}{ccl}
S_n \odot S_m & = & \sum_{k \geq 0} a^{nk} \odot \sum_{k' \geq 0}
a^{mk'} = \sum_{p\geq 0} \langle S_n|a^p\rangle \langle S_m|a^p\rangle a^p \\
& = & \sum_{k \geq 0}
a^{lcm(n,m)k}= S_{lcm(n,m)}\enspace .
\end{array}
$$
Thus, for $n$ and $m$ coprime, the rank of the product is $nm$, which proves the claim. \qed

\section{Shuffle of automata compatible with relators}
In this section, we deal with automata whose actions can be coded by
elements of a
monoid defined by
generators and relations. The first interesting case historically
encountered is the trace monoid but, as we
will see below, some results can be extended to the general case.
To end with, we study the relators permitting the shuffle of automata.
\subsection{Series over a monoid and automata}
In the whole section $R\subset A^*\times A^*$ is a relator and $\equiv_R$
is the congruence relation generated by $R$.
\begin{defn}\label{DMini1}
\begin{enumerate}
\item
Let $f:A^* \rightarrow X$ ($X$ a set) and $\equiv$ be a congruence on $A^*$, we will say that
$f$ is $\equiv-
$compatible if
$$
u\equiv v\Rightarrow f(u)=f(v).
$$
\item An automaton ${\cal A}=(\lambda,\mu,\gamma)$ is said
$\equiv$-compatible if
$\mu:A^*\rightarrow K^{n\times n}$ is.
\end{enumerate}
\end{defn}
\begin{remarks}\label{RMini1}
\begin{enumerate}
\item The coarsest congruence compatible with a function $f$ is known as the
syntactic congruence of $f$. A
 non trivial result says that the syntactic congruence
of all Greene's invariants is the plactic equivalence \cite{SchD}.\\
\item \label{pt2} If an automaton ${\cal A}$ is $\equiv$-compatible, then it is
straigthforward to see that its behaviour is.\\
\item We can restate geometrically (\ref{pt2}) of definition \ref{DMini1} as
:
\begin{center}
For each state $q$ and $(u,v)\in R$ then $q.u=q.v$.
\end{center}
\item \label{pt4}If $f:A^*\rightarrow M$ is a morphism of monoids ( this is the
case for the data  $\mu$ of automata ) compatibility has just to be tested on $R$, more precisely
$$
(\forall(u,v)\in R)(f(u)=f(v))\Rightarrow f
\mbox{ is }\equiv\mbox{-compatible}.
$$
\item If $S,T$ are $\equiv$-compatible, so is $S\odot T$ (which is by no means the case for $\shuffl$ and $\uparrow$, see
discussion below).
\end{enumerate}
\end{remarks}
The converse of remark \ref{RMini1}(\ref{pt2}) is true for minimal automata over
fields as shown just below.
\begin{prop}\label{PComp}
Suppose that $K$ is a field (commutative or skew). \\
Let $S: A^*\rightarrow K$ be a rational series,
the following assertions are equivalent:
\begin{enumerate}
\item \label{16pt1}$S$ is $\equiv$-compatible.\\
\item \label{16pt2}The minimal automata of $S$ are $\equiv$-compatible.
\end{enumerate}
\end{prop}
{\bf Proof} Let us first prove that (\ref{16pt1})$\Rightarrow$(
\ref{16pt2}). By the minimality of ${\cal A}$, it exists words $u_1,
u_2,\dots, u_n,v_1, v_2,\dots v_n$ such that the column block matrix
$L=(\lambda\mu(u_i))_{i\in[1,n]}$ and the line block matrix
$R=(\mu(v_i)\gamma)_{i\in[1,n]}$ are invertible $n\times n$ matrices
($K$ may not be commutative see \cite{Flo}). Thus, if $w\equiv w'$
then
$$
\begin{array}{rcl}
L\mu(w)R&=&(\lambda\mu(u_iwv_j)\gamma)_{1\leq i,j\leq n}\\
&=&(\langle S|u_iwv_j\rangle)_{1\leq i,j\leq n}\\
&=&(\langle S|u_iw'v_j\rangle)_{1\leq i,j\leq n}\\
&=&(\lambda\mu(u_iw'v_j)\gamma)_{1\leq i,j\leq n}\\
&=&L\mu(w')R
\end{array}
$$
And thus, $\mu(w)=\mu(w')$.\\
The converse is straightforward from remark \ref{RMini1}(\ref{pt4}).\qed

It is clear that $\equiv$-compatibility is stable under linear combinations
(i.e. if the series $(S_{i,j})_{(i,j)\in I\times J}$ are $\equiv$-compatible so is
$\sum\alpha_iS_{i,j}\beta_j$). However, the Cauchy product of two compatible
series may not be so, as shown by the  example:
$ab\equiv ba$, $S=a$ and $T=b$.\\
\subsection{Study for general semirings}
In case of a field, the compatibility of automata with shuffle product is
equivalent to the compatibility of the coproduct with the congruence and its
square. More precisely
\begin{thm}
\begin{enumerate}
\item \label{17pt1}
Suppose that $K$ is a field. Let $\equiv$ be a congruence with finite fibers\footnote{i.e. the
classes of $\equiv$ are finite sets.},
the following assertions are equivalent.
\begin{enumerate}
\item \label{17pta}If ${\cal A}_1$ and ${\cal A}_2$ are two $\equiv$-compatible
automata so is ${\cal A}_1\mbox{ }\scriptshuffl \mbox{ }{\cal A}_2$.
\item \label{17ptb}The coproduct respects $\equiv$ in the following sense:\\
For every $(u,v)\in A^*\times A^*$, we have
$$
u\equiv v\Rightarrow c(u)\equiv^{\otimes 2} c(v).
$$
\end{enumerate}
where $\equiv^{\otimes2}$ is the "square" of $\equiv$ defined as the kernel
of the natural mapping
$$
K \langle A \rangle \otimes K \langle A \rangle \rightarrow K[A^*/_\equiv]\otimes K[A^*/_\equiv].
$$
\item \label{17pt2}The preceding conditions imply that if $S$ and $T$ are two
$\equiv$-compatible series, so are $S\shuffl T$, $S\uparrow T$.
\end{enumerate}
\end{thm}
{\bf Proof} To prove (\ref{17pt1}.\ref{17ptb}) $\Rightarrow$
(\ref{17pt1}.\ref{17pta}), it suffices to remark that
$\mu=(\mu_1\otimes \mu_2)\circ c$ where $\mu_1,\mu_2$ and $\mu$ are
respectively the associated morphisms of the automata ${\cal A}_1$,
${\cal A}_2$ and ${\cal A}_1\mbox{ }\scriptshuffl\mbox{ }{\cal
A}_2$.\\
Now, we prove that (\ref{17pt1}.\ref{17pta}) $\Rightarrow$ (\ref{17pt1}.\ref{17ptb}).
We consider the (product order) relation on the multidegrees ( $\alpha,\beta\in\bbbn^{(A)}$ ):
$$
(\alpha\leq\beta)\Leftrightarrow (\forall
a\in A)(\alpha(a)\leq\beta(a)).
$$
Let $w$ be a word. In the sequel, we denote $[w]$ the mapping $(a\rightarrow |w|_a)$ its 
multidegree and
$Cl(w)$ its equivalence class modulo $\equiv$. Let $w_1\equiv w_2$ be two
equivalent words.
Consider
$$
t_1=\sup_{w\in Cl(w_1)}[w].
$$
And let ${\cal C}_1\dots{\cal
C}_k$ be
the classes which contain at least a word whose multidegree is less than
$t_1$, and we set
$$
t_2=\sup_{w\in\cup_{i=1}^k{\cal C}_i}[w]
$$
($t_1$ and $t_2$ are well defined due to the "finite fibers" hypothesis).\\
With $A^{\leq t_2}:=\{w/[w]\leq t_2\}$, let us define the following truncation of $\equiv$ by
$$
u\sim v\Leftrightarrow\left\{
\begin{array}{c}
Cl(u)\not\subseteq A^{\leq t_2}\mbox{ and } Cl(v)\not\subseteq A^{\leq t_2}\\
\mbox{ or }\\
Cl(u)=Cl(v)
\end{array}
\right.
$$
The following lemma is easy.
\begin{lem}\label{LTh1|1}
\begin{enumerate}
\item The equivalence $\sim$ is a congruence coarser than $\equiv$.
\item The classes of $\sim$ are ${\cal C}_1,{\cal C}_2,\dots,{\cal C}_k,{\cal
C}_{k+1},\dots {\cal C}_{p-1}
$ and $${\cal C}_p=\bigcup_{Cl(w)\not\subseteq A^{\leq t_2}}Cl(w)$$ where
${\cal C}_1,\dots,{\cal C}_{p-1}$
are equivalence classes of $\equiv$ precisely the equivalence classes of $\equiv$ which are 
subsets of
$A^{\leq t_2}$.
\item In particular $w_1\sim w_2$ and $[w_i]\leq t_1$ implies $w_1\equiv w_2$.
\end{enumerate}
\end{lem}
For every $a\in A$, we define $\mu(a)$ as the matrix (with respect to the
basis $({\cal C}_j)_{j\in[1,p]}$) of the linear
transformation $\OV u\rightarrow \OV a.\OV u\in A^*/_\sim$, where $\OV u$ denotes the class of $u$ for $\sim$. More
explicitly
$$
\mu(w): {\cal C}_j\rightarrow \OV w.{\cal C}_j.
$$
Then, $\mu$ is $\equiv$-compatible and hence the automata ${\cal
A}_{i,j}=(e_{{\cal C}_i},\mu,e^*_{{\cal C}_j})$ (with $(e_{{\cal
C}_{i}})_{1\leq i\leq p}$ being the
canonical basis of $K^{p\times 1}$) are $\equiv$-compatible. Then, by
(\ref{17pta}) the $p^4$ automata
$$
{\cal A}_{i_1,j_1}\mbox{ }\scriptshuffl\mbox{ }{\cal
A}_{i_2,j_2}=(e_{{\cal C}_{i_1}}\otimes e_{{\cal C}_{i_2}},\mu\otimes
I_{p}+I_{p}\otimes\mu,e^*_{{\cal C}_{j_1}}\otimes e^*_{{\cal C}_{j_2}})
$$
are $\sim$-compatible. This, implies that the morphism $\nu:A^*\rightarrow
K^{p^2\times p^2}$ defined by $\nu(a)=\mu(a)\otimes
I_{p}+I_{p}\otimes \mu(a)$ for each $a\in A$, is $\sim$-compatible. Now, as
$w_1\equiv w_2$, one has
$$
\begin{array}{rcl}
\sum_{I+J=[1\dots n]}\mu(w_1[I])\otimes\mu(w_1[J])&=&\nu(w_1)\\
&=&\nu(w_2)\\
&=&\sum_{I+J=[1\dots n]}\mu(w_2[I])\otimes\mu(w_2[J])
\end{array}
$$
which proves (evaluating this linear transformation on $1\otimes 1$) that
$$
\sum_{I+J=[1\dots n]}w_1[I]\otimes w_1[J]\sim^{\otimes2}\sum_{I+J=[1\dots
n]}w_2[I]\otimes w_2[J]
$$
but, as $[w_i[I]],[w_i[J]]\leq t_1$ for $I,J \subset [1..n]$, lemma \ref{LTh1|1} implies
$c(w_1)\equiv^{\otimes2} c(w_2)$.\\ \\
Now, we prove (\ref{17pt1}) $\Rightarrow$ (\ref{17pt2}). In fact we have, $\langle S\shuffl
T|w\rangle =\langle S\otimes T|c(w)\rangle$. As $S$ and $T$ are $\equiv$-compatible, the
assertion (\ref{17pt1}.\ref{17ptb}) implies the $\equiv$-compatibility of $S\shuffl T$.\qed

In fact (\ref{17pt1}.\ref{17ptb}) can be formulated without the hypothesis over $K$ and the
fibers of $\equiv$ and then (\ref{17pt1}.\ref{17ptb}) $\Rightarrow$ (\ref{17pt1}.\ref{17pta}) in the (very) general
case.\\
According to this remark we can give the following definition.
\begin{defn}
Let $K$ be a semiring. A congruence will be said $K-\shuffl$ compatible if
(\ref{17pt1}.\ref{17ptb}) is fullfilled.
\end{defn}
Partial commutations are $K-\shuffl$ compatible for any $K$, so does
, more generally, the relators $a^{p^{e_1}}b^{p^{e_2}}\equiv b^{p^{e_2}}a^{p^{e_1}}$ and
$a^{p^{e_1}}=b^{p^{e_2}}$
for $K=\bbbz/p\bbbz$ with $p$ prime.\\ \\
In the next paragraph we completely solve the problem of $K-\shuffl$
compatibility for semirings which are not rings.\\
The case when $K$ is a ring of
characteristic $0$ is known (see \cite{DK1}) but the tools developped below
shows this again by a different argument.
\subsection{Generalities}
In the following we need some elementary properties.
\begin{lem}\label{Lgeneral1}
Let $\phi: K_1\rightarrow K_2$ be a morphism of semirings then
\begin{enumerate}
\item If $\equiv$ is $K_1-\shuffl$ compatible then it is $K_2-\shuffl$
compatible.
\item If $\phi$ is into, the converse is true.
\end{enumerate}
\end{lem}
{\bf Proof} Straightforward, remarking that the mapping
$\bbbn.1_{K_1}\displaystyle\mathop\rightarrow^\phi \bbbn.1_{K_2}$ is
surjective.\qed

\begin{rem}
This lemma implies that if a congruence is $\bbbn-\shuffl$ compatible then
it
is $K-\shuffl$
compatible for each semiring $K$. In fact, a congruence is
$K-\shuffl$ compatible if and only if it is $\bbbn.1_K-\shuffl$ compatible.
\end{rem}
Let $K$ be a semiring, in the following we discuss according to the subsemiring
$K_0=\bbbn.1_K$. The semiring $K_0$ is entirely characterized by the monoid structure of
$(K_0,+)$ which depends of the two following parameters:
$$
m(K)=inf\{e\in \bbbn/\exists r\in \bbbn^*,e.1_k=(e+r).1_K\}\in\bbbn\cup\{+\infty\}
$$
and if $m(K)\neq\infty$
$$
l(K)=inf\{r\in \bbbn^*/m(K).1_K=(m(K)+r).1_K\}\in \bbbn^*.
$$
\begin{lem}\label{Lgeneral2}
Let $R$ be a relator on $A^*$. Then, $\equiv_R$ is $K-\shuffl$ compatible
if and only
if for each pair $(w_1,w_2)\in R$ we have $c(w_1)\equiv^{\otimes2}_Rc(w_2)$.
\end{lem}
{\bf Proof} The "if" part is straightforward considering the
morphism
$$
c:A^*/_\equiv\rightarrow K[A^*/_\equiv]\otimes K[A^*/_\equiv].
$$
The converse is obvious.\qed
\begin{lem}\label{Lgeneral4}
Each congruence generated by relators under the form $a\equiv b$ or
$cd\equiv dc$ with
$a,b,c,d\in A$ is $K-\shuffl$ compatible.
\end{lem}
{\bf Proof} According to lemma \ref{Lgeneral2}, it suffices to check
that
$$
c(a)=a\otimes1+1\otimes a\equiv^{\otimes2} b\otimes 1+1\otimes b=c(b)
$$
for each $a\equiv b\in A$ and
$$
\begin{array}{ll}
c(cd) & =cd\otimes1+c\otimes d+d\otimes c+1\otimes cd \\
& \equiv^{\otimes2}
dc\otimes1+c\otimes d+d\otimes c+1\otimes dc\\
& =c(dc)
\end{array}
$$
for each pair of letters $(a,b)\in A^2$ such that $cd\equiv dc$.\qed

\begin{lem}\label{Lgeneral5}
Let $B\subseteq A$ be a subalphabet. If $\equiv$ is $K-\shuffl$ compatible
then so is the congruence
$\equiv_B:=\equiv\cap B^2$ .
\end{lem}
{\bf Proof} Direct computation.\qed\\ \\ The following general lemma
will be used later.
\begin{lem}\label{Lgeneral3}
Let $u\in A^+$ be a word and let $n$ be the maximal integer such that $u$
can be
written under the form
$u=u_1a^n$ with $u_1\in A^*$, $a\in A$ and $n\geq 1$ then
$$
\langle c(u)|u_1\otimes a^n\rangle=1.
$$
\end{lem}
{\bf Proof} Suppose that $n=1$ then it is easy to verify that
$u_1\otimes a$ appears only one times in the polynomial $c(u)$. By
induction on $n$, we find the result.\qed

\subsection{The case when $m(K)\neq 0$}
{ a) The boolean case}

We first consider the case where $K=\bbbb$ is the boolean semiring. The
$\bbbb-\shuffl$ compatible congruences are caracterised by the following
result.
\begin{prop}\label{PBool1}
A congruence is $\bbbb-\shuffl$ compatible if and only if it is generated by
the following
relators
$$
\left\{\begin{array}{lc}
a\equiv 1&(LE)\\
a\equiv b&(LI)\\
ab\equiv ba&(LC)
\end{array}\right.
$$
\end{prop}
{\bf Proof} Let us first prove that a congruence is $\bbbb-\shuffl$
compatible if it is generated by relators (LE), (LI) or (LC).
According lemmas \ref{Lgeneral2} and \ref{Lgeneral4}, it suffices to
prove that the relators (LE) are $\bbbb-\shuffl$ compatible. In
fact, we have
$$
a\equiv 1\Rightarrow c(a)=a\otimes 1+1\otimes a\equiv^{\otimes2}
(1+1)\otimes
1=c(1)
$$
which proves the result.\\
Now, we prove the converse. Let $A'=\{a\in A/a\not\equiv1\}$ and $S\subseteq A'$ be a section of
$\equiv\cap A'\times A'$. It is clear that if (LE) is a list of couples $\{(a,1)\}_{a\in A-A'}$
and (LI) a list of couples $\{(a,b)\}_{x\equiv y, x\in S, y\in A'-S}$, then $\equiv$ is generated
by $\equiv_S:=\equiv\cap S^*\times S^*$, (LI) and (LE). So, it suffices to prove that $\equiv_S$
is generated by $(LC)$ relators. Let us prove first, that $\equiv_S$ is multihomogeneous. Let
$\equiv_m$ be the
multihomogeneous part of $\equiv_S$ (i.e. the congruence generated by the
pairs
$(u,v)\in\equiv_S$ such that $[u]=[v]$). Let $(u,v)$ be a pair of
words such that $u\equiv_S v$ and $u\not\equiv_m v$ with $|u|$ minimal.
Suppose that $u=1$, if $v\neq 1$ we can set $v=v_1a$ with $a\in S$. Then, as by lemma
\ref{Lgeneral5} $\equiv_S$ is again $\bbbb-\shuffl$ compatible,
$$
\langle\OV v_1\otimes \OV a|c(\OV1)\rangle=1
$$
( $\OV w$ denoting the class of $w$ for $\equiv$), but $c(1)=1\otimes1$ which implies $a\equiv_S 1$ and contradicts the
construction of $S$. Then, $u\neq 1$ and we can write $u$ under the form
$u=u_1a$ with $a\in S$. As $\langle c(u)|u_1\otimes a\rangle=1$, it exists two
complementary subwords $v[I]$ and $v[J]$ of $v$ such that $v[I]\otimes
v[J]\equiv_S^{\otimes2} u_1\otimes a$. But, $v\equiv_Su_1a\equiv_Sv[I]v[J]$ which implies
$v\equiv_mv[I]v[J]$ and proves $\equiv_S=\equiv_m$.\\
Let $\equiv_\theta$ be the congruence generated by pairs $(ab,ba)$ with
$a,b\in S$ and $ab\equiv_S ba$.

\begin{lem}\label{LBool2}
Let $u\equiv_S v$ with $v\in S^*a$ then it exists $u_1\equiv_\theta u$ with
$u_1\in S^*a$.
\end{lem}
{\bf Proof} We have $[u]=[v]$ from what precedes and in particular
$|u|_a\neq0$. Let $u_1=u_2au'_2$ be a word such that
$u_1\equiv_\theta u$, $|u'_2|_a=0$ and $|u'_2|$ minimal. Suppose
that $u'_2\neq 1$, then we can write $u'_2=bu_3$ with $b\in S$ and
$u_3\in S^*$. Let $a^qb=(u_1)_I$ be the subword of $u_1$ with $q$
maximal ($q=|u_a|$, the word is unique but the equality has
$|u'_2|_b$ solutions in I), it exists two complementary subwords
$v[I]$ and $v[J]$ such that $a^qb\otimes w\equiv_S^{\otimes2}
v[I]\otimes v[J]$ where $w$ is a subword of $u$ complementary of
$a^qb$. Then $a^qb\equiv_S v[I]$ and then, as
$|v[I]|_a=|u|_a=|v|_a$, $v[I]=a^{q-i}ba^i $ with $i\geq1$. This
implies $ab\otimes a^{q-1}\equiv_m^{\otimes2} ab\otimes
a^{q-1}+ba\otimes a^{q-1}$. As $\equiv_S$ is multihomogeneous, we
have necessarily $ab\equiv_S ba$. It follows $u\equiv_\theta
u_2abu_3\equiv_\theta u_2bau_3$ which contradicts the minimality of
$|u'_2|$ and proves the result.\qed\\ \\
{\bf End of the proof of proposition} If $\equiv_S\not=
\equiv_\theta$, let $(u,v)$ be a couple of words such that
$u\equiv_Sv$ and $u\not\equiv_\theta v$ with $|u|+|v|$ minimal.\\
Let $a$ be a letter such that
$u\equiv_\theta u_1a^k=u', v\equiv_\theta v_1a^l=v'$ with
$k,l\not=0, k+l\geq 2$ maximal (the existence of a such letter follows from
lemma \ref{LBool2}). Without restriction
we can suppose that $k\leq l$. We have $\langle u_1\otimes a^k|c(u')\rangle=1$ and then
it exists two complementary subwords $v'[I]$ and $v'[J]$ of $v'$ such that
$u_1\otimes a^k\equiv^{\otimes 2}_Sv'[I]\otimes v'[J]$. Hence, the
multihomogeneity of $\equiv_S$ gives $v'[J]=a^k$ and we can write
$v'[I]=v_2a^\alpha$ where $v_2$ is a subword of $v_1$. If $\alpha>0$, we
have $u_1\equiv_S v_2a^\alpha$ and by lemma \ref{LBool2}, it would exist
$u_2\in S^*$ such that $u_1\equiv_\theta u_2a$. Hence,
$u\equiv_\theta u_2a^{k+1}$ which contradicts the maximality of $k+l$.
Thus $\alpha=0$ and $v'[I]\notin S^*a$ is a subword of $v_1$, we have thus
$|u|-k=|u_1|=|v'[I]|\leq |v_1|=|v|-l$ but we had $k\leq l$ then $k=l$. Now
$v_1=v'[I]$ and then
$u_1\equiv_\theta v_1$ which implies
$$
u\equiv_\theta u_1a^k\equiv_\theta v_1a^k\equiv_\theta v
$$
a contradiction, this proves the result.\qed\\ \\
 { b) Other
semirings such that $m(K)\neq 0$}

\begin{thm}\label{TOther}
Let $K$ be a semiring such that $m(K)\neq 0$. Then a congruence $\equiv$ is
$K-\shuffl$ compatible if
and only if
\begin{enumerate}
\item \label{28pti} If $1_K+1_K=1_K$, it is generated by relators (LE), (LI) and (LC).
\item \label{28ptii} If $1_K+1_K\neq 1_K$, it is generated by relators (LI) and (LC).
\end{enumerate}
In the two cases, $A^*/\equiv$ is a partially commutative monoid.
\end{thm}
{\bf Proof} The assertion (\ref{28pti}) can be easily proved using
lemma \ref{Lgeneral1} and proposition \ref{PBool1}. Let us show the
assertion (\ref{28ptii}). Let $K$ be a semiring such that $m(K)\neq
0$ and $1_K+1_K\neq 1_K$, then it exists a morphism from $K$ onto
$\bbbb$ (this morphism sends $0$ on $0$ and $x\neq 0$ on $1$). Let
$\equiv$ be a $K-\shuffl$ compatible congruence, by lemma
\ref{Lgeneral1} $\equiv$ is so $\bbbb-\shuffl$ compatible and then
it is generated by (LE), (LI) or (LC) relators. A fast computation
shown that the only possibilities are (LI) and (LC). Which gives the
result.\qed
\begin{cor}\cite{DK1}
Let $K$ be a ring of characteristic $0$. A congruence is $K-\shuffl$ compatible if and only if it 
is generated by relators of the type (LI) and (LC).
\end{cor}
\begin{example}
Let $\bbbn_{max}=(\bbbn\cup\{-\infty\},max,+)$ be the tropical semiring and
$A=\{a,b,c,d\}$, the congruence generated by $\{(a,1),(a,b),(cd,dc)\}$ is
$\bbbn_{ max}-\shuffl$ compatible.
\end{example}
{ c) Other examples in characteristic $2$}

We consider here the field $K=\bbbz/2\bbbz$, and the relators $$
R=\{(ab^2,b^2a),(a^2b,ba^2),(abab,baba)\}.
$$
It is obvious to
see that the congruence generated by the set $\{(ab^2,b^2a)$, $(a^2b,ba^2)\}$ is $\bbbz/2\bbbz-\shuffl$ compatible. 
Furthermore,
we have
$$
\begin{array}{rcl}
c(abab)&=&abab\otimes 1+aba\otimes b+a^2b\otimes b+ab^2\otimes a+bab\otimes a+ba\otimes ab\\
& + & va^2\otimes b^2 + ab\otimes ba+b^2\otimes a^2+b\otimes aba+b\otimes a^2b+a\otimes b^2a\\
& + & a\otimes bab+1\otimes abab \\
&\equiv_R^{\otimes2}&baba\otimes 1+aba\otimes b+ba^2\otimes b+b^2a\otimes a+bab\otimes a+ba\otimes ab\\
& + & a^2\otimes b^2 + ab\otimes ba+b^2\otimes a^2+b\otimes aba+b\otimes ba^2+a\otimes ab^2\\
& + & a\otimes bab+1\otimes baba\\
&=&c(baba)
\end{array}
$$
which implies the $\bbbz/2\bbbz-\shuffl$ compatibility of $\equiv_R$.\\
We can remark that this property does not occur if $K$ is not a field or if $2_K\neq 0_K$.\\
In the same way, the congruence generated by the relators
$$
R=\{(a^8b^2,b^2a^8),(a^4b^4,b^4a^4),(a^4b^2a^4b^2,b^2a^4b^2a^4)\}
$$
is $\bbbz/2\bbbz-\shuffl$ compatible.
\section{Conclusion}
Many computations over rational series can be lifted at the level of automata
and these (classical) constructions has been proved to be genericaly optimal.
The implementation of classical rational laws ( shuffle, Hadamard, infiltration) has
suggested us other laws (which also preserve rationality) and we have proved
that, under some natural hypothesis, there is no other choice than a deformation of the classical
case.\\
The study of the shuffle product over automata raises the question of the
compatibility with relators. The answer is of course coefficient dependant and in
classical cases ($0$ characteristic, boolean and proper semirings) it is
interesting to observe that only dependance relations can occur. But the
$p$-characteristic induces strange phenomena and opens some new and exciting questions.

%


\end{document}